\documentclass[11pt,leqeq]{article}
\usepackage{amssymb,amsthm}
\usepackage{amsmath}
\usepackage{amscd}
\usepackage{xy}
\xyoption{all}
\newtheorem{thm}{Theorem}
\newtheorem{lema}[thm]{Lemma}
\newtheorem{prop}[thm]{Proposition} 

\newtheorem{defi}[thm]{Definition}

\begin{document}
\title{q,k-generalized gamma and beta functions}
\author{Rafael D\'\i az\thanks{Work partially supported by UCV.}\ \
and Carolina Teruel\thanks{Work partially supported by IVIC.}}
\maketitle
\begin{abstract}
We introduce the q,k-generalized Pochhammer symbol. We construct
$\Gamma_{q,k}$ and $B_{q,k}$, the q,k-generalized
 gamma and beta fuctions, and show that they satisfy properties that generalize those satisfied by the classical gamma and beta functions.
Moreover, we provide integral representations for $\Gamma_{q,k}$
and $B_{q,k}.$
\end{abstract}

\section{Introduction}
The usefulness of the gamma and beta functions can hardly be
overstated. However, recent results coming from the combinatorics
of annihilation and creation operators \cite{QSF}, and the
construction of hypergeometric functions from the point of view of
the k-generalized Pochhammer symbol given in \cite{ED}, have impel
upon us the need to introduce subtle, yet deep, generalizations of
the all mighty gamma and beta functions. The main goal of this
paper is to introduce a two parameter deformation of the classical
gamma and beta functions, which we call the q,k-generalized gamma
and beta functions and will be denoted by $\Gamma_{q,k}$ and
$B_{q,k}$, respectively. $\Gamma_{q,k}$ and $B_{q,k}$ fit into the
following conmutative diagrams

$$\xymatrix{ \Gamma_{q,k}(t) \ar[d]_{k \rightarrow 1} \ar[r]^{q
\rightarrow 1}
& \Gamma_{k}(t)  \ar[d]^{k \rightarrow 1} \\
\Gamma_{q}(t) \ar[r]_{q \rightarrow 1} & \Gamma(t)  } \quad \quad
\quad  \quad \quad \quad \xymatrix{ B_{q,k}(t,s) \ar[d]_{k
\rightarrow 1} \ar[r]^{q \rightarrow 1}
& B_{k}(t,s)  \ar[d]^{k \rightarrow 1} \\
B_{q}(t,s) \ar[r]_{q \rightarrow 1} & B(t,s)  }$$

Let us explain the notation used in the diagrams above.  Recall
that the Euler's gamma and beta functions are given by the
following Riemann integrals:

$$\Gamma(t)=\int_{0}^{\infty} x^{t-1} e^{-x} dx,\quad t>0.$$

$$B(t,s)=\int_{0}^{1} x^{t-1}(1-x)^{s-1}dx = \int_{0}^{\infty}
\frac{x^{t-1}}{(1+x)^{t+s}} dx, \quad t,s>0.$$

Our motivation to introduce the q,k-generalized gamma and beta
functions is two-fold: on one hand, D\'\i az and Pariguan in
\cite{ED} defined  a k-deformation, $k>0$ a real number, of the
gamma and beta functions given by the following Riemann integrals

\begin{equation*}
\Gamma_k(t)=\int_{0}^{\infty} x^{t-1} e^{-\frac{x^k}{k}}dx, \quad
t>0.
\end{equation*}
\begin{equation*}
B_k(t,s)= \int_{0}^{\infty} x^{t-1}(1+x^k)^{-\frac{t+s}{k}}dx,
\quad t,s>0.
\end{equation*}

On the other hand, De Sole and Kac in \cite{So} introduced a
q-deformation, $0<q<1$  a real number, of the gamma and beta
functions given by the following Jackson integrals

\begin{equation} \label{q-gamma}
\Gamma_q(t)=\int_{0}^{\infty/(1-q)} x^{t-1} E^{-qx} d_qx, \quad
t>0.
\end{equation}
\begin{equation}\label{q-beta}
B_q(t,s)=\int_{0}^{1} x^{t-1}(1-qx)_q^{s-1}d_qx, \quad t,s>0.
\end{equation}

Although equation (\ref{q-beta}) has been known for about a
century, the upper limit in formula (\ref{q-gamma}) has only been
recently established by Koelink and Koornwinder \cite{K},
\cite{koor}: the factor $\frac{1}{(1-q)}$ was traditionally
omitted yielding  a divergent Jackson integral.

Our two-parameter deformations $\Gamma_{q,k}$ and $B_{q,k}$
generalize both constructions above. We show that our function
$\Gamma_{q,k}$ is related to the q,k-generalized Pochhammer symbol
$[t]_{n,k}$, to be defined in Section \ref{Basic}, in the same way
as the classical $\Gamma$ function and the Pochhammer symbol are
related to each other. Our functions $\Gamma_{q,k}$ and $B_{q,k}$
are given by the following formulae

\begin{equation*}
\displaystyle{\Gamma_{q,k}(t)=\frac{{(1-q^k)_{q,k}^{{\frac{t}{k}}
-1}}}{{(1-q)^{\frac{t}{k}-1}}}}, \quad t>0.
\end{equation*}

\begin{equation*}
\displaystyle{B_{q,k}(t,s)=\frac{(1-q)(1-q^k)_{q,k}^{\frac{s}{k}-1}}{(1-q^t)_{q,k}^{\frac{s}{k}}}},\quad
\mbox{ for all } s,t>0.
\end{equation*}

In Section 4, we give Jackson integral representations for the
q,k-generalized gamma and beta functions in terms of the
q,k-analogue of the exponential function $E_{q,k}^x$. These are
given by

\begin{equation*}
\displaystyle{\Gamma_{q,k}(t)=\int_{0}^{\left(\frac{[k]_q}{(1-q^k)}\right)^{\frac{1}{k}}}x^{t-1}E_{q,k}^{-\frac{q^k
x^k}{[k]_q}}d_qx}, \quad t>0.
\end{equation*}

\begin{equation*}
\displaystyle{B_{q,k}(t,s)=[k]_q^{-\frac{t}{k}}\int_{0}^{[k]_{q}^{\frac{1}{k}}}x^{t-1}\left(
1-q^k\frac{x^k}{[k]_q}\right)_{q,k}^{\frac{s}{k}-1}d_qx,}\quad
t,s>0.
\end{equation*}

Furthermore, in Section 5, we give other integral representations
for $\Gamma_{q,k}$ and $B_{q,k}$, using the q,k-analogue of the
exponential function $e_{q,k}^x$. These integral representations
are given by

\begin{equation}\label{alanis}
\displaystyle{\Gamma_{q,k}(t)=c(a,t)\int_{0}^{\infty/a(1-q^k)^{\frac{1}{k}}}
x^{t-1} e_{q,k}^{-\frac{x^k}{[k]_q}}}d_qx, \quad t>0.
\end{equation}

\begin{equation}\label{morissette}
\displaystyle{B_{q,k}(t,s)=c(a,t)[k]_q^{-\frac{t}{k}}
\int_{0}^{\infty /a} \frac{x^{t-1}}{\left(1+
\frac{x^k}{[k]_q}\right)_{q,k}^{\frac{t+s}{k}}} d_qx, \quad
t,s>0,}
\end{equation}
where
$\displaystyle{c(a,t)=\frac{a^t[k]_q^{\frac{t}{k}}}{1+[k]_qa^{k}}\left(
1+ \frac{1}{[k]_qa^k}\right)_{q,k}^{\frac{t}{k}} \left( 1+[k]_qa^k
\right)_{q,k}^{1-\frac{t}{k}}}.$

\vspace{0.5cm} We shall see that formulae (\ref{alanis}) and
(\ref{morissette}) are deeply related to the classical Jacobi
triple product identity and to the famous Ramanujan identity,
respectively. We remark that all our results are elementary and
require little knowledge beyond a  basic introduction to
q-calculus.

\section{Basic results}\label{Basic}

In this Section we introduce the q,k-generalized Pochhammer symbol
and also some basic definitions that will be used in the rest of
the paper. For completeness we review well-known material that
might be found, for example, in \cite{Ch}, \cite{Pa}, \cite{HY}
and \cite{RA}. Let us begin by introducing q-derivatives and
Jackson integrals, see \cite{To},\cite{Jac}.

\begin{defi}
\newcommand{\Func}{\operatorname{Func}}
Let us denote by $\Func(\mathbb{R},\mathbb{R})$ the real vector
space of all functions from $\mathbb{R}$ to $\mathbb{R}$. Fix
$0<q<1$ and consider the linear operators
$I_q,d_q,\partial_q:\Func(\mathbb{R},\mathbb{R})\longrightarrow
\Func(\mathbb{R},\mathbb{R})$, given for all
$f\in\Func(\mathbb{R},\mathbb{R})$ by:

\begin{itemize}
\item $I_q(f)(x)=f(qx)$, \quad for all $x\in\mathbb{R}$.\\
\item $d_q(f)= I_q(f)-f$.\\
\item $\displaystyle{\partial_q(f)=\frac{d_qf}{d_qx}=\frac{I_q(f)-
f}{(q-1)x}}$.\ \ \ \   $\partial_q(f)$ is called the q-derivative
of the
function f.\\
\end{itemize}

\end{defi}
\begin{defi}
\newcommand{\Func}{\operatorname{Func}}
\begin{enumerate}
\item The definite q-integral of a function $f \in
\Func(\mathbb{R},\mathbb{R})$ from $0$ to $b>0$ is given by
$$\int_{0}^{b}f(x)d_qx=(1-q)b\sum_{n=0}^{\infty}q^nf(q^nb).$$\\
\vspace{-0.5cm} \item The improper q-integral of a function $f \in
\Func(\mathbb{R},\mathbb{R})$ is given by

$$ \int_0^{\infty/a} f(x) d_qx =
(1-q)\sum_{n\in\mathbb{Z}}\frac{q^n}{a}f\left(\frac{q^n}{a}\right).$$\\
\end{enumerate}
\end{defi}
\vspace{-1cm}

\begin{prop}
\newcommand{\Func}{\operatorname{Func}}
For any functions $ f,g \in \Func(\mathbb{R},\mathbb{R})$, the
following properties hold:
\begin{itemize}
\item $\partial_q (f+g)= \partial_q(f) + \partial_q(g).$\\
\item{$\partial_q(fg)= f \partial_q (g) + I_q (g)\partial_q (f)$.}\\
\item$\displaystyle \partial_q (f/g)= \frac{\partial_q (f)g
- f\partial_q(g)}{I_q (g)(g)}.$\\
\item $\partial_q(f(a x^{b}))= a[b]_q(\partial_{q^{b}}f)(a
x^{b})x^{b-1}, \quad \mbox{ for all } a,b \in \mathbb{R}$.\\
\item $f(b)g(b)-f(a)g(a)= \displaystyle{\int_a^b f
\partial_q(g)d_qx + \int_a^b
I_q(g)\partial_q(f)d_qx},$ \quad for all \ $0\leq a < b \leq +\infty$.\\
\end{itemize}
\end{prop}

\begin{defi}
Let $0<q<1$ be a fixed real number. Let us denote by ${[\mbox{
}]_q: \mathbb{R} \rightarrow \mathbb{R}}$ the map given by
$\displaystyle{[t]_q=\frac{(1-q^t)}{(1-q)}},$ for all $t \in
\mathbb{R}.$
\end{defi}

The map ${[\mbox{ }]_q: \mathbb{R} \rightarrow \mathbb{R}}$ it is
not an algebra homomorphism. Nevertheless, the following
identities are satisfied:
\begin{enumerate}
\begin{minipage}[v]{10cm}
  \item $\displaystyle{[s+t]_q=[s]_q + q^s[t]_q,}$ for all $s,t \in \mathbb{R}.$ \item $\displaystyle{[st]_q=[s]_{q^t} [t]_q,}$ for all $s,t \in \mathbb{R}.$
\end{minipage}
\begin{minipage}[v]{10cm}
  \item $\displaystyle{[1]_q=1.}$\item $\displaystyle{[0]_q=0.}$
\end{minipage}

\end{enumerate}

Next definition is fundamental for the rest of the paper. Indeed,
our original motivation for this work was to find integral
representations for the q,k-generalized Pochhammer symbol.
\begin{defi}

\begin{enumerate}
    \item Let $t \in \mathbb{R}$ and $n \in \mathbb{Z^+}$. The
k-generalized Pochhammer symbol is given by \vspace{-0.5cm}
$$(t)_{n,k}= t(t+k)(t+2k). \ . \ .(t+(n-1)k)=\prod_{j=0}^{n-1}(t+jk).$$\\
\vspace{-1.3cm} \item The q,k-generalized Pochhammer symbol is
given by \vspace{-0.5cm}
$$[t]_{n,k}= [t]_q[t+k]_q[t+2k]_q. \ . \ .[t+(n-1)k]_q =
\prod_{j=0}^{n-1}[t+jk]_q.$$\\
\end{enumerate}
\end{defi}
\vspace{-1.3cm} Notice that $[t]_{n,k}\rightarrow(t)_{n,k}$ as
$q\rightarrow1.$ Let us introduce some notation that will be used
throughout the paper.

\begin{defi}\label{sepa}
Let $x,y,t\in\mathbb{R}$ and $n \in \mathbb{Z^+}$
\begin{enumerate}
    \item $\displaystyle{(x+y)_{q,k}^n := \prod_{j=0}^{n-1} (x+ q^{jk}y).}$\\
    \item
$\displaystyle{(1+x)_{q,k}^{t}:=\frac{(1+x)_{q,k}^{\infty}}{(1+q^{kt}x)_{q,k}^{\infty}}.}$
\end{enumerate}
\end{defi}
\begin{lema}\label{separa}
Let $x,s,t\in\mathbb{R}$. Then
$(1+x)_{q,k}^{s+t}=(1+x)_{q,k}^{s}(1+q^{ks}x)_{q,k}^{t}.$
\end{lema}

\section{Explicit formulae for $\Gamma_{q,k}$ and
$B_{q,k}$} The k-generalization of the gamma function introduced
by D\'\i az and Pariguan in \cite{ED}, is univocally determined by
the following properties: \vspace{0.3cm}
\begin{enumerate}
\begin{minipage}[c]{6cm}
\item $\Gamma_k(t+k)=t\Gamma_k(t), \quad t>0.$
\end{minipage}
\begin{minipage}[c]{3cm}
\item $\Gamma_k(k)=1.$
\end{minipage}
\begin{minipage}[t]{10cm}
\item $\Gamma_k $ is logarithmically convex.\\
\end{minipage}
\end{enumerate}
Properties 1 and 2 imply that $\displaystyle{(t)_{n,k}=
\frac{\Gamma_k(t+nk)}{\Gamma_k(t)}}$, \ \ for all $t>0$ and
$n\in\mathbb{Z^+}.$

We define the q,k-generalized gamma function $\Gamma_{q,k}$ by
demanding it satisfies the q,k-analogue of properties 1 and 2
above. Thus we assume that $\Gamma_{q,k}$ is such that:
$\Gamma_{q,k}(t+k)=[t]_q \Gamma_{q,k}(t)$ and $\Gamma_{q,k}(k)=1.$
This implies that,
\begin{eqnarray*}
\Gamma_{q,k}(nk)=\displaystyle{\prod_{j=1}^{n-1}[jk]_q}=
\displaystyle{ \prod_{j=1}^{n-1} \frac{(1-q^{jk})}{{(1-q)}}}
=\displaystyle{\frac{(1-q^k)_{q,k}^{n-1}}{(1-q)^{n-1}}}.
\end{eqnarray*}
After the change of variable t=nk, one is lead to the following
\begin{defi}\label{defigamma}
The function $\Gamma_{q,k}$ is given by the formula

\begin{equation*}
\displaystyle{\Gamma_{q,k}(t)=\frac{{(1-q^k)_{q,k}^{{\frac{t}{k}}
-1}}}{{(1-q)^{\frac{t}{k}-1}}}}, \quad t>0.
\end{equation*}
\end{defi}

\begin{lema}\label{ginfini}
The infinite product expression for the function $\Gamma_{q,k}$ is
given by
\begin{equation*}
\Gamma_{q,k}(t)=\frac{(1-q^k)_{q,k}^{\infty}}{(1-q^t)_{q,k}^{\infty}(1-q)^{\frac{t}{k}-1}},
\quad t>0.
\end{equation*}
\end{lema}

Next proposition guarantees that $\Gamma_{q,k}$ indeed satisfies
the q,k-analogue of properties 1 and 2 above.

\begin{prop}\label{propiedades}
The function $\Gamma_{q,k}$ satisfies the following identities for
$t>0$:
\begin{enumerate}
    \item \label{1}$\Gamma_{q,k}(t+k)=[t]_q \Gamma_{q,k}(t).$\\
    \item \label{2}$\Gamma_{q,k}(k)=1.$\\
    \item \label{3}$\displaystyle{ \frac{\Gamma_{q,k}(t+nk)}{\Gamma_{q,k}(t)}= [t]_{n,k}},$\quad for all $n\in\mathbb{Z^+}.$ \\

\end{enumerate}

\begin{proof}
\ref{1}.
$\displaystyle{\Gamma_{q,k}(t+k)=\frac{(1-q^t)}{(1-q)}\frac{(1-q^k)_{q,k}^{\frac{t}{k}-1}}{(1-q)^{\frac{t}{k}-1}}=[t]_q\Gamma_{q,k}(t).}$
\vspace{0.5cm}

\ref{2}. Obvious. \vspace{0.5cm}

\ref{3}.
$\displaystyle{\frac{\Gamma_{q,k}(t+nk)}{\Gamma_{q,k}(t)}=
\frac{(1-q^t)_{q,k}^{n}}{(1-q)^n} =
\prod_{j=0}^{n-1}\frac{(1-q^{t+jk})}{(1-q)}=
\prod_{j=0}^{n-1}[t+jk]_q = [t]_{n,k}.}$
\end{proof}
\end{prop}

\begin{defi}\label{defi}
The function $B_{q,k}(t,s)$ is given by the formula

\begin{equation*}\label{q,k-beta}
\displaystyle{B_{q,k}(t,s)=\frac{\Gamma_{q,k}(t)\Gamma_{q,k}(s)}{\Gamma_{q,k}(t+s)},
\quad \mbox{ for all } s,t>0.}
\end{equation*}
\end{defi}
Which in turns imply the next
\begin{lema}\label{betica}
\begin{enumerate}
\item
$\displaystyle{B_{q,k}(t,s)=\frac{(1-q)(1-q^k)_{q,k}^{\frac{s}{k}-1}}{(1-q^t)_{q,k}^{\frac{s}{k}}}},\quad
\mbox{ for all } s,t>0.$\\
\item
$\displaystyle{B_{q,k}(t,s)=\frac{(1-q)(1-q^k)_{q,k}^{\infty}(1-q^{s+t})_{q,k}^{\infty}}{(1-q^s)_{q,k}^{\infty}(1-q^t)_{q,k}^{\infty}}},
\quad \mbox{ for all } s,t>0.$\\

\end{enumerate}
\begin{proof}
Use Definition \ref{defigamma} and Definition \ref{defi}.
\end{proof}
\end{lema}

\begin{prop}\label{propo}
The function $B_{q,k}$ satisfies the following formulae for
$s,t>0$
\begin{enumerate}
    \item \label{4}$\displaystyle{B_{q,k}(t,\infty)=(1-q)^{\frac{t}{k}}\Gamma_{q,k}(t)}.$ \\
    \item \label{5}$\displaystyle{B_{q,k}(t+k,s)= \frac{[t]_q}{[s]_q} B_{q,k}(t,s+k)}.$ \\
    \item \label{6}$\displaystyle{B_{q,k}(t,s+k)=B_{q,k}(t,s) - q^s B_{q,k}(t+k,s) }.$\\
    \item \label{7}$\displaystyle{B_{q,k}(t,s+k)=\frac{[s]_q}{[s+t]_q}B_{q,k}(t,s)}.$\\
    \item \label{8}$\displaystyle{B_{q,k}(t,k)= \frac{1}{[t]_q}}.$\\
    \item \label{9}$\displaystyle{B_{q,k}(t,nk)=
    (1-q)\frac{(1-q^k)_{q,k}^{n-1}}{(1-q^t)_{q,k}^{n}}=
    (1-q)
    \frac{(1-q^k)_{q,k}^{n-1}(1-q^k)_{q,k}^{\frac{t}{k}-1}}{(1-q^k)_{q,k}^{\frac{t}{k}+n-1}}},\quad n\in\mathbb{Z^+}.$\\
\end{enumerate}

\begin{proof}
\begin{enumerate}
\item
$\displaystyle{B_{q,k}(t,\infty)=\frac{(1-q)(1-q^k)_{q,k}^{\infty}}{(1-q^t)_{q,k}^{\infty}}=\frac{(1-q)(1-q^k)_{q,k}^{\infty}}{(1-q^{k(\frac{t-k}{k})}q^k)_{q,k}^{\infty}}}
=(1-q)(1-q^k)_{q,k}^{\frac{t-k}{k}}.$\\
Then $\displaystyle{B_{q,k}(t,\infty)={(1-q)^{\frac{t}{k}}}\Gamma_{q,k}(t)}.$\\
\vspace{0.2cm}

\item
$\displaystyle{\frac{B_{q,k}(t+k,s)}{B_{q,k}(t,s+k)}=\frac{(1-q^k)_{q,k}^{\frac{s}{k}-1}(1-q^t)_{q,k}^{\frac{s}{k}+1}}{(1-q^k)_{q,k}^{\frac{s}{k}}(1-q^{t+k})_{q,k}^{\frac{s}{k}}}=\frac{(1-q^t)}{(1-q^s)}=\frac{[t]_q}{[s]_q}}.$\\

\vspace{0.5cm} \item $\:$ \vspace{-1.2cm}
\begin{eqnarray*}
\frac{B_{q,k}(t,s+k)-B_{q,k}(t,s)}{B_{q,k}(t+k,s)}&=&\frac{(1-q^{t+k})_{q,k}^{\frac{s}{k}}(1-q^k)_{q,k}^{\frac{s}{k}}}{(1-q^t)_{q,k}^{\frac{s}{k}+1}(1-q^k)_{q,k}^{\frac{s}{k}-1}}-\frac{(1-q^{t+k})_{q,k}^{\frac{s}{k}}}{(1-q^t)_{q,k}^{\frac{s}{k}}}\\
&=&\frac{(1-q^s)}{(1-q^t)}-\frac{(1-q^{t+s})}{(1-q^t)}=-q^s.\hspace{6cm}
\end{eqnarray*}\\

 \item $\:$ \vspace{-1.2cm}
\begin{eqnarray*}\frac{B_{q,k}(t,s+k)}{B_{q,k}(t,s)}&=&
\frac{(1-q^k)_{q,k}^{\frac{s}{k}}(1-q^t)_{q,k}^{\frac{s}{k}}}{(1-q^t)_{q,k}^{\frac{s}{k}+1}(1-q^k)_{q,k}^{\frac{s}{k}-1}}\\
&=&\frac{(1-q^k)_{q,k}^{\frac{s}{k}}(1-q^t)_{q,k}^{\frac{s}{k}}}{(1-q^t)_{q,k}^{\frac{s}{k}}(1-q^{t+s})(1-q^k)_{q,k}^{\frac{s}{k}}(1-q^s)^{-1}}\\
&=&\frac{(1-q^s)}{(1-q^{t+s})}=\frac{[t]_q}{[t+s]_q}.\hspace{8cm}
\end{eqnarray*}\\

\vspace{-0.6cm}
\item$\displaystyle{B_{q,k}(t,k)=\frac{(1-q)}{(1-q^t)}=\frac{1}{[t]_q}}.$\\

 \item Using  items \ref{7} and \ref{8} of this
proposition, we obtain:
\begin{eqnarray*}
B_{q,k}(t,nk)=
\displaystyle{\frac{\displaystyle{\prod_{j=1}^{n-1}}[jk]_q}{\displaystyle{\prod_{j=0}^{n-1}}[jk+t]_q}}
= \frac{(1-q)(1-q^k)_{q,k}^{n-1}}{(1-q^t)_{q,k}^{n}}=
(1-q)\frac{(1-q^k)_{q,k}^{n-1}(1-q^k)_{q,k}^{\frac{t}{k}-1}}{(1-q^k)_{q,k}^{\frac{t}{k}+n-1}}.\\
\end{eqnarray*}
Letting $n \rightarrow \infty$, we get
 $B_{q,k}(t,\infty)=(1-q)(1-q^k)_{q,k}^{\frac{t}{k}-1}.$\\
\end{enumerate}
\end{proof}
\end{prop}

\section{Integral representations for $\Gamma_{q,k}$ and
$B_{q,k}$}\label{rep} As promised in the introduction, in this
Section we provide Jackson integral representations for our
$\Gamma_{q,k}$ and $B_{q,k}$ functions, in terms of the
q,k-analogue
exponential function $E_{q,k}^{x}$. Recall that \\
$$\displaystyle{E_{q,k}^{x}=\sum_{n=0}^{\infty}\frac{q^{kn(n-1)/2}x^n}{[n]_{q^k}!}=(1+(1-q^k)x)_{q,k}^{\infty}}.$$

\begin{thm}\label{tildes}
\begin{enumerate}
\item The function $\Gamma_{q,k}(t)$ is given by the following
Jackson integral
\begin{equation}\label{gamatilde}
\displaystyle{\Gamma_{q,k}(t)=\int_{0}^{\left(\frac{[k]_q}{(1-q^k)}\right)^{\frac{1}{k}}}x^{t-1}E_{q,k}^{-\frac{q^k
x^k}{[k]_q}}d_qx}, \quad t>0.
\end{equation}\\

\vspace{-0.7cm} \item The function $B_{q,k}$ is given by the
following Jackson integral
\begin{equation}\label{betatilde}
\displaystyle{B_{q,k}(t,s)=[k]_q^{-\frac{t}{k}}\int_{0}^{[k]_{q}^{\frac{1}{k}}}x^{t-1}\left(
1-q^k\frac{x^k}{[k]_q}\right)_{q,k}^{\frac{s}{k}-1}d_qx}, \quad
t,s>0.
\end{equation}\\
\end{enumerate}
\end{thm}
\vspace{-0.8cm} In order to prove Theorem \ref{tildes} we begin by
denoting the right hand side of (\ref{gamatilde}) by
$\overline{\Gamma}_{q,k}$, and the right hand side of
(\ref{betatilde}) by $\overline{B}_{q,k}$. Let us check that
$\overline{\Gamma}_{q,k}$ and $\overline{B}_{q,k}$ satisfy
properties analogue  to those stated in Proposition
\ref{propiedades} and Proposition \ref{propo} for $\Gamma_{q,k}$
and $B_{q,k}$, respectively.

\begin{prop}\label{propis}
The function $\overline{\Gamma}_{q,k}$ satisfies the following
identities for $t>0$
\begin{enumerate}
\item \label{2}$\overline{\Gamma}_{q,k}(k)=1.$\\
\item \label{1}$\overline{\Gamma}_{q,k}(t+k)=[t]_q \overline{\Gamma}_{q,k}(t).$\\
\end{enumerate}

\begin{proof}

\begin{enumerate}

\item $\displaystyle{\overline{\Gamma}_{q,k}(k)=
-\int_{0}^{\left(\frac{[k]_q}{(1-q^k)}\right)^{\frac{1}{k}}}
\partial_q\left(
E_{q,k}^{\frac{x^k}{[k]_q}} \right)d_qx=1,}$ since $E_{q,k}^{-\frac{1}{1-q^k}}=0,$ and $E_{q,k}^0=1.$\\

\item  Using q-integration by parts
\begin{eqnarray*}
\overline{\Gamma}_{q,k}(t+k)&=&\int_{0}^{\left(\frac{[k]_q}{(1-q^k)}\right)^{\frac{1}{k}}}x^{t+k-1}E_{q,k}^{\frac{-q^kx^k}{[k]_q}}d_qx\\
&=&
-\int_{0}^{\left(\frac{[k]_q}{(1-q^k)}\right)^{\frac{1}{k}}}x^{t}\partial_q
\left(E_{q,k}^{-\frac{x^k}{[k]_q}}\right)d_qx\\
&=&[t]_q\int_{0}^{\left(\frac{[k]_q}{(1-q^k)}\right)^{\frac{1}{k}}}x^{t-1}E_{q,k}^{-\frac{q^kx^k}{[k]_q}}d_qx.
\hspace{5cm}\\\nonumber
\end{eqnarray*}
\end{enumerate}
\vspace{-0.6cm}
\end{proof}
\end{prop}
\vspace{-0.5cm}
\begin{prop}\label{propos}
The function $\overline{B}_{q,k}$ satisfies the following formulae
for $s,t>0$
\begin{enumerate}
    \item \label{4}$\displaystyle{\overline{B}_{q,k}(t,\infty)=(1-q)^{\frac{t}{k}}\overline{\Gamma}_{q,k}(t)}.$ \\
    \item \label{5}$\displaystyle{\overline{B}_{q,k}(t+k,s)= \frac{[t]_q}{[s]_q} \overline{B}_{q,k}(t,s+k)}.$ \\
    \item \label{6}$\displaystyle{\overline{B}_{q,k}(t,s+k)=\overline{B}_{q,k}(t,s) - q^s \overline{B}_{q,k}(t+k,s) }.$\\
    \item \label{7}$\displaystyle{\overline{B}_{q,k}(t,s+k)=\frac{[s]_q}{[s+t]_q}\overline{B}_{q,k}(t,s)}.$\\
    \item \label{8}$\displaystyle{\overline{B}_{q,k}(t,k)= \frac{1}{[t]_q}}.$\\
     \item \label{9}$\displaystyle{\overline{B}_{q,k}(t,nk)=
    (1-q)\frac{(1-q^k)_{q,k}^{n-1}}{(1-q^t)_{q,k}^{n}}=
    (1-q)
    \frac{(1-q^k)_{q,k}^{n-1}(1-q^k)_{q,k}^{\frac{t}{k}-1}}{(1-q^k)_{q,k}^{\frac{t}{k}+n-1}}},\quad n\in\mathbb{Z^+}.$\\
\end{enumerate}

\begin{proof}
\begin{enumerate}
\item Using the change $x=(1-q^k)^{\frac{1}{k}}y$, (\ref{ppp}) is
obtained from (\ref{pp}). \vspace{-0.3cm}
\begin{eqnarray}\label{pp}
\overline{B}_{q,k}(t,\infty)&=&[k]_q^{-\frac{t}{k}}\int_{0}^{[k]_{q}^{\frac{1}{k}}}x^{t-1}E_{q,k}^{-\frac{q^kx^k}{(1-q^k)[k]_q}}d_qx\\\label{ppp}
&=&(1-q)^{\frac{t}{k}}\int_{0}^{\left(\frac{[k]_q}{(1-q^k)}\right)^{\frac{1}{k}}}y^{t-1}E_{q,k}^{-\frac{q^k
y^k}{[k]_q}}d_qy\\
&=&(1-q)^{\frac{t}{k}}\overline{\Gamma}_{q,k}(t).\nonumber
\end{eqnarray}\\

\vspace{-0.6cm}
 \item Using the formula
$\partial_q\left(1+b\frac{x^k}{[k]_q}\right)_{q,k}^{t}=\frac{[kt]_q}{[k]_q}bx^{k-1}\left(1+bq^k\frac{x^k}{[k]_q}\right)_{q,k}^{t-1}$
in (\ref{co}) we have (\ref{ck})
\begin{eqnarray}\label{co}
 \overline{B}_{q,k}(t+k,s)&=&[k]_{q}^{-\frac{t}{k}-1}\int_{0}^{[k]_{q}^{\frac{1}{k}}}x^{t+k-1}\left(1-q^k\frac{x^k}{[k]_q}\right)_{q,k}^{\frac{s}{k}-1}d_qx\\\label{ck}
&=&-\frac{[k]_q^{-\frac{t}{k}}}{[s]_q}\int_{0}^{[k]_{q}^{\frac{1}{k}}}x^t
\partial_q\left(\left(1-\frac{x^k}{[k]_q}\right)_{q,k}^{\frac{s}{k}}
\right)d_qx\\
&=&\frac{[t]_q}{[s]_q}\overline{B}_{q,k}(t,s+k).\nonumber
\end{eqnarray}\\
\vspace{-0.6cm}
\item Using Lemma \ref{separa}, we get
\begin{align*}
\overline{B}_{q,k}(t,s+k)&=[k]_q^{-\frac{t}{k}}\int_{0}^{[k]_{q}^{\frac{1}{k}}}x^{t-1}\left(
1-q^k\frac{x^k}{[k]_q}\right)_{q,k}^{\frac{s}{k}-1}\\\nonumber
&-q^s[k]_q^{-\frac{t}{k}-1}\int_{0}^{[k]_{q}^{\frac{1}{k}}}x^{t+k-1}\left(
1-q^k\frac{x^k}{[k]_q}\right)_{q,k}^{\frac{s}{k}-1}d_qx\\
\nonumber
&=\overline{B}_{q,k}(t,s)-q^s\overline{B}_{q,k}(t+k,s).\nonumber
\end{align*}\\
\vspace{-0.6cm}
\item It is easy to check that
$\overline{B}_{q,k}(t,s+k)=\frac{[s]_q}{[s+t]_q}\overline{B}_{q,k}(t,s)$ using the properties 2 and 3 above.\\

\item
$\displaystyle{\overline{B}_{q,k}(t,k)=[k]_q^{-\frac{t}{k}}\int_{0}^{[k]_{q}^{\frac{1}{k}}}x^{t-1}=\frac{1}{[t]_q}.}$\\

\item Use items 4 and 5 of this proposition.\\
\end{enumerate}\end{proof}

\end{prop}
\vspace{-0.4cm} \textit{Proof of Theorem \ref{tildes}.} By
Proposition \ref{propo} part 1,
$\displaystyle{\overline{B}_{q,k}(t,\infty)=(1-q)^{\frac{t}{k}}\overline{\Gamma}_{q,k}(t)},
\  t>0.$ Also, it follows from Proposition \ref{propo} part 6 that
$\displaystyle{\overline{B}_{q,k}(t,\infty)=(1-q)^{\frac{t}{k}}\Gamma_{q,k}(t)},
\  t>0.$ Therefore, $$\overline{\Gamma}_{q,k}(t)=\Gamma_{q,k}(t),
\ t>0.$$

Items 3 and 4 of Proposition \ref{propos} imply that
$$\overline{B}_{q,k}(t,s)=B_{q,k}(t,s), \mbox{ for all } t>0 \mbox{ and } s=nk,
 \quad \mbox{with } n\in\mathbb{Z^+}.$$ We would like to prove that
\begin{equation*}\label{q,k-beta}
\displaystyle{\overline{B}_{q,k}(t,s)=\frac{\Gamma_{q,k}(t)\Gamma_{q,k}(s)}{\Gamma_{q,k}(t+s)},
\quad \mbox{ for all } s,t>0.}
\end{equation*}
By Lemma \ref{ginfini} we have that
\begin{equation*}
\frac{\Gamma_{q,k}(t)\Gamma_{q,k}(s)}{\Gamma_{q,k}(t+s)}=\frac{(1-q)(1-q^k)_{q,k}^{\infty}(1-q^{t+s})_{q,k}^{\infty}}{(1-q^t)_{q,k}^{\infty}(1-q^s)_{q,k}^{\infty}}.
\end{equation*}
\vspace{0.4cm}

By the definite Jackson integral and Definition \ref{sepa} item 2,
we obtain \vspace{-0.3cm}
\begin{eqnarray*}\nonumber
\overline{B}_{q,k}(t,s)&=&[k]_q^{-\frac{t}{k}}\int_{0}^{[k]_{q}^{\frac{1}{k}}}x^{t-1}\left(
1-q^k\frac{x^k}{[k]_q}\right)_{q,k}^{\frac{s}{k}-1}d_qx\\\nonumber
&=&(1-q)\sum_{n=0}^{\infty}q^{nt}(1-q^{k(n+1)})_{q,k}^{\frac{s}{k}-1}\\
&=&(1-q)\sum_{n=0}^{\infty}q^{nt}
\frac{(1-q^{k(n+1)})_{q,k}^{\infty}}{(1-q^{s+nk})_{q,k}^{\infty}}.
\end{eqnarray*}
So, we have to show that
$$\frac{(1-q)(1-q^k)_{q,k}^{\infty}(1-q^{t+s})_{q,k}^{\infty}}{(1-q^t)_{q,k}^{\infty}(1-q^s)_{q,k}^{\infty}}=(1-q)\sum_{n=0}^{\infty}q^{nt}
\frac{(1-q^{k(n+1)})_{q,k}^{\infty}}{(1-q^{s+nk})_{q,k}^{\infty}}.$$

Making the changes $u=q^t$ and $v=q^s$, we reduce our problem to
prove that
\begin{equation}\label{sides}
\frac{(1-q)(1-q^k)_{q,k}^{\infty}(1-uv)_{q,k}^{\infty}}{(1-u)_{q,k}^{\infty}(1-v)_{q,k}^{\infty}}=(1-q)\sum_{n=0}^{\infty}b^n
\frac{(1-q^{k(n+1)})_{q,k}^{\infty}}{(1-vq^{nk})_{q,k}^{\infty}}.
\end{equation}

Now, both sides of equation (\ref{sides}) are formal power series
in q with rational coefficients in $u$ and $v$. Since we already
know that they agree for an infinite number of values, namely
$u=q^t$ and $v=q^s$, where $t>0$ and $s=nk$ with
$n\in\mathbb{Z^+}$, the desired result holds.

\section{Other integral representations}
In this Section we provided Jackson integral representations for
the $\Gamma_{q,k}$ and $B_{q,k}$ using the q,k-analogue of the
exponential function $e_{q,k}^x$. We remark that in an earlier
version of this paper we took the upper limit in our integrals to
be $\infty$ leading to divergent integrals. The correct upper
limit in our definition was clear to us only after reading
\cite{So}. Recall that
$$\displaystyle{e_{q,k}^x=\sum_{n=0}^{\infty}\frac{x^n}{[n]_{q^k}!}=\frac{1}{(1-(1-q^k)x)_{q,k}^{\infty}}}.$$

\begin{defi}
The function $\gamma_{q,k}^{(a)}$, $a>0,$ is given by the
following Jackson integral
$$\displaystyle{\gamma_{q,k}^{(a)}(t)=\int_{0}^{\infty/a(1-q^k)^{\frac{1}{k}}} x^{t-1} e_{q,k}^{-\frac{x^k}{[k]_q}}}d_qx, \quad t>0.$$
\end{defi}

Next proposition shows that $\gamma_{q,k}^{(a)}$ satisfies
properties similar to those given in Proposition \ref{propiedades}
for the $\Gamma_{q,k}$ function.
\begin{prop}
The function $\gamma_{q,k}^{(a)}$ satisfies the following formulae
for $a,t>0$
\begin{enumerate}
    \item \label{10}$\gamma_{q,k}^{(a)}(k)=1.$\\
    \item \label{11}$ \displaystyle{\gamma_{q,k}^{(a)}(t+k)=q^{-t} [t]_q \gamma_{q,k}^{(a)}(t).}$\\
    \item \label{12}$\displaystyle{\gamma_{q,k}^{(a)}(nk)=q^{-kn(n-1)/2}\Gamma_{q,k}(nk)}$, \quad for every $n \in \mathbb{Z^+}.$\\
\end{enumerate}

\begin{proof}
\begin{enumerate}
\item
$\displaystyle{\gamma_{q,k}^{(a)}(k)=\int_{0}^{\infty/a(1-q^k)^{\frac{1}{k}}}
x^{k-1}
e_{q,k}^{-\frac{x^k}{[k]_q}}d_qx=-\int_{0}^{\infty/a(1-q^k)^{\frac{1}{k}}}\partial_q\left(
 e_{q,k}^{-\frac{x^k}{[k]_q}}\right)=1.}$\\

\vspace{0.5cm}

\item $\:$
\vspace{-1.2cm}
\begin{eqnarray*}
\gamma_{q,k}^{(a)}(t+k)&=&-q^{-t}\int_{0}^{\infty/a(1-q^k)^{\frac{1}{k}}}(qx)^t
\partial_q\left(e_{q,k}^{-\frac{x^k}{k}}\right)d_qx\\
&=&[t]_qq^{-t}\int_{0}^{\infty/a(1-q^k)^{\frac{1}{k}}} x^{t-1}
e_{q,k}^{-\frac{x^k}{[k]_q}}d_qx =
q^{-t}[t]_q\gamma_{q,k}^{(a)}(t). \hspace{18cm}
\end{eqnarray*}\\
\vspace{-0.7cm}
\item From items \ref{10} and \ref{11} above, we have that\\
$$\displaystyle{\gamma_{q,k}^{(a)}(nk)=q^{-kn(n-1)/2}\prod_{j=1}^{n-1}[jk]_q=q^{-kn(n-1)/2}\Gamma_{q,k}(nk).}$$\\
\end{enumerate}
\vspace{-0.6cm}
\end{proof}
\end{prop}
\vspace{-0.5cm}
\begin{defi}
The function $\beta_{q,k}^{(a)}$, $a>0,$ is given by the following
Jackson integral
$$\displaystyle{\beta_{q,k}^{(a)}(t,s)=[k]_q^{-\frac{t}{k}} \int_{0}^{\infty /a}
\frac{x^{t-1}}{\left(1+
\frac{x^k}{[k]_q}\right)_{q,k}^{\frac{t+s}{k}}} d_qx, \quad
t,s>0.}$$
\end{defi}

\begin{prop}\label{propo2}
The function $\beta_{q,k}^{(a)}$ satisfies the following formulae
for $a,t,s>0$
\begin{enumerate}
    \item \label{13}$\displaystyle{\beta_{q,k}^{(a)}(t,\infty)=(1-q)^{\frac{t}{k}}\gamma_{q,k}^{(a)}(t).}$\\
    \item \label{14}$\displaystyle{\beta_{q,k}^{(a)}(t+k,s)= q^{-t}\frac{[t]_q}{[t+s]_q}\beta_{q,k}^{(a)}(t,s) }$.\\
    \item \label{15}$\displaystyle{\beta_{q,k}^{(a)}(k,s)= \frac{1}{[s]_q}.} $\\
    \item \label{16}$ \displaystyle{\beta_{q,k}^{(a)}(t,s+k)=\frac{[s]_q}{[t+s]_q} \beta_{q,k}^{(a)}(t,s)}.$\\
    \item \label{17}$\beta_{q,k}^{(a)}(nk,s)=q^{-kn(n-1)/2}B_{q,k}(nk,s)$, for all $n \in \mathbb{Z^+}.$\\
\end{enumerate}

\begin{proof}
\begin{enumerate}
\item Using the change $x=(1-q^k)^{\frac{1}{k}}y$,  (\ref{ecc}) is
obtained from (\ref{ec}).

\begin{eqnarray}\label{ec}
\beta_{q,k}^{(a)}(t,\infty)&=&[k]_q^{-\frac{t}{k}}
\int_{0}^{\infty /a}x^{t-1}
e_{q,k}^{-\frac{x^k}{(1-q^k)[k]_q}}d_qx. \\\label{ecc}
&=&[k]_q^{-\frac{t}{k}}(1-q^k)^{\frac{t}{k}}
\int_{0}^{\infty/a(1-q^k)^{\frac{1}{k}}}y^{t-1}e_{q,k}^{-\frac{y^k}{[k]_q}}d_qy\\\label{eccc}
&=&(1-q)^{\frac{t}{k}}\gamma_{q,k}^{(a)}(t).\nonumber
\end{eqnarray}
\item Using the formula $\partial_q\left( \frac{\left(
1+a\frac{x^k}{[k]_q}\right)_{q,k}^{s}}{\left(
1+b\frac{x^k}{[k]_q}\right)_{q,k}^{t}}\right)=
    \frac{[ks]_qax^{k-1}\left(1+aq^k\frac{x^k}{[k]_q}\right)_{q,k}^{s-1}}{[k]_q\left(
    1+bq^k\frac{x^k}{[k]_q}\right)_{q,k}^{t}}
    -\frac{[kt]_qbx^{k-1}\left( 1+a\frac{x^k}{[k]_q}\right)_{q,k}^{s}}{[k]_q \left(1+b\frac{x^k}{[k]_q}
    \right)_{q,k}^{t+1}}$.
\begin{eqnarray*}
\beta_{q,k}^{(a)}(t+k,s)&=&-\frac{[k]_q^{-\frac{t}{k}}}{[t+s]_q}q^{-t}\int_{0}^{\infty/a}(qx)^{t}\partial_q\left(\frac{1}{\left(
1+\frac{x^k}{[k]_q} \right)_{q,k}^{\frac{t+s}{k}}}\right)d_qx\\
&=&\frac{[t]_q}{[t+s]_q}q^{-t}[k]_q^{-\frac{t}{k}}\int_{0}^{\infty/a}\frac{x^{t-1}}{\left(
1+\frac{x^k}{[k]_q}\right)_{q,k}^{\frac{t+s}{k}}}d_qx\\
&=&\frac{[t]_q}{[t+s]_q}q^{-t}\beta_{q,k}^{(a)}(t,s).\\
\end{eqnarray*}
\item Using the identity $\partial_q\left(
\frac{1}{\left(1+\frac{x^k}{[k]_q}
\right)_{q,k}^{\frac{s}{k}}}\right)=-\frac{x^{k-1}[s]_q}{[k]_q\left(1+\frac{x^k}{[k]_q}
\right)_{q,k}^{\frac{s}{k}+1}}$, we get
\begin{eqnarray*}
\beta_{q,k}^{(a)}(k,s)=-\frac{1}{[s]_q}\int_{0}^{\infty
/a}\partial_q\left( \frac{1}{\left( 1+\frac{x^k}{[k]_q}
\right)_{q,k}^{\frac{s}{k}}}\right)d_qx=\frac{1}{[s]_q}.
\end{eqnarray*}\\
\vspace{-0.6cm}
\item Using $\partial_q\left(
\frac{ax^{ks}}{[k]_q^{s}\left(1+b\frac{x^k}{[k]_q}\right)_{q,k}^{t}}\right)=\frac{ax^{ks-1}[ks]_q}{[k]_q^{s}\left(1+b\frac{x^k}{[k]_q}\right)_{q,k}^{t+1}}
- b([kt]_q-[ks]_q)
\frac{ax^{k(s+1)-1}}{[k]_q^{s+1}\left(1+b\frac{x^k}{[k]_q}\right)_{q,k}^{t+1}}.$\\
\begin{eqnarray*}
\beta_{q,k}^{(a)}(t,s+k)&=&\frac{[k]_q^{\frac{t+s}{k}}[k]_q^{-\frac{t}{k}}q^s}{[t+s]_q}\int_{0}^{\infty
/a} \frac{1}{(qx)^{s}} \partial_q\left(
\frac{x^{t+s}}{[k]_q^{\frac{t+s}{k}}\left(
1+\frac{x^k}{[k]_q}\right)_{q,k}^{\frac{t+s}{k}}}\right)d_qx\\
&=&-\frac{q^s[k]_q^{-\frac{t}{k}}}{[t+s]_q}\int_{0}^{\infty
/a}\frac{x^{t+s}}{\left(1+
\frac{x^k}{[k]_q}\right)_{q,k}^{\frac{t+s}{k}}}\partial_q\left(
\frac{1}{x^s}
\right)d_qx\\
&=&\frac{[s]_q}{[t+s]_q}[k]_q^{-\frac{t}{k}}\int_{0}^{\infty
/a}\frac{x^{t-1}}{(1+ \frac{x^k}{[k]_q})_{q,k}^{\frac{t+s}{k}}}
d_qx=\frac{[s]_q}{[t+s]_q}\beta_{q,k}^{(a)}(t,s).
\end{eqnarray*}\\
\vspace{-0.8cm}
\item Using property 2 above recursively\\
$\displaystyle{\beta_{q,k}^{(a)}(nk,s)=q^{-kn(n-1)/2}\frac{(1-q)(1-q^k)_{q,k}^{n-1}(1-q^{nk})_{q,k}^{\frac{s}{k}-n}}{(1-q^s)_{q,k}^{n}(1-q^{nk})_{q,k}^{\frac{s}{k}-n}}=q^{-kn(n-1)/2}B_{q,k}(nk,s)}.$\\
\end{enumerate}
\end{proof}
\end{prop}

\begin{lema}\label{lemota}
Let $s,t \in \mathbb{R}$ and $n\in\mathbb{Z^+}$, we have the
following identities
\begin{enumerate}
\item $\displaystyle{(1+q^{ks}x)_{q,k}^{t}=\frac{(1+x)_{q,k}^{s+t}}{(1+x)_{q,k}^{s}}=\frac{(1+x)_{q,k}^{t}(1+q^{kt}x)_{q,k}^{s}}{(1+x)_{q,k}^{s}}.}$\\
\item $\displaystyle{(1+q^{-kn}x)_{q,k}^{t}=(1+x)_{q,k}^{t}\frac{(x+q^k)_{q,k}^{n}}{(q^{kt}x+q^k)_{q,k}^{n}}.}$\\
\end{enumerate}
\end{lema}
Next theorem provides our second integral representation for the
functions $\Gamma_{q,k}$ and $B_{q,k}$.

\begin{thm}\label{teor}
For every $a,s,t>0$ we have:
\begin{enumerate}
\item $\Gamma_{q,k}(t)=c(a,t)\gamma_{q,k}^{(a)}(t).$\\
\item$B_{q,k}(t,s)=c(a,t)\beta_{q,k}^{(a)}(t,s).$\\
\end{enumerate} Where
$$\displaystyle{c(a,t)=\frac{a^t[k]_q^{\frac{t}{k}}}{1+[k]_qa^{k}}\left(
1+ \frac{1}{[k]_qa^k}\right)_{q,k}^{\frac{t}{k}} \left( 1+[k]_qa^k
\right)_{q,k}^{1-\frac{t}{k}}}.$$

\begin{proof}
Since both $\displaystyle{B_{q,k}(t,s+k)=\frac{[s]_q}{[s+t]_q}
B_{q,k}(t,s)}$ and $
\displaystyle{\beta_{q,k}^{(a)}(t,s+k)=\frac{[s]_q}{[t+s]_q}
\beta_{q,k}^{(a)}(t,s).}$ It is clear that if $c(a,t)$ is such
that
$c(a,t)\beta_{q,k}^{(a)}(t,s)=B_{q,k}(t,s)$, then c(a,t) must be\\
$$\left(\displaystyle{\int_{0}^{\infty/a}\partial_q\left(
\frac{x^t}{[k]_q{^\frac{t}{k}}\left(
1+\frac{x^k}{[k]_q}\right)_{q,k}^{\frac{t}{k}}}
\right)d_qx}\right)^{-1}.$$

We know that
$\displaystyle{\beta_{q,k}^{(a)}(t,k)=\frac{1}{[t]_q}\displaystyle{\int_{0}^{\infty/a}\partial_q\left(
\frac{x^t}{[k]_q{^\frac{t}{k}}\left(
1+\frac{x^k}{[k]_q}\right)_{q,k}^{\frac{t}{k}}} \right)d_qx}}$.
Thus by definition of q-derivative and the Jackson integral, we
have:

\begin{equation}\label{c}
\displaystyle{\int_{0}^{\infty/a}\partial_q (F_k) d_qx=
\lim_{n\rightarrow \infty}F_k \left(\frac{1}{aq^n}
\right)-\lim_{n\rightarrow \infty}F_k\left(\frac{q^n}{a} \right),}
\end{equation}

\vspace{0.5cm} where the limits are taken over integers.

From (\ref{c}), we obtain:
$\displaystyle{\beta_{q,k}^{(a)}(t,k)=\frac{1}{[t]_q}\lim_{n\rightarrow
\infty}\left(
[k]_q^{\frac{t}{k}}(aq^n)^t\left(1+\frac{1}{[k]_q(aq^n)^k}\right)_{q,k}^{\frac{t}{k}}\right)^{-1}.}$

Using Lemma \ref{lemota} part 2 in equation (\ref{hola}) we have
\begin{eqnarray}\label{hola}
c(a,t)&=&[k]_q^{\frac{t}{k}}a^t\displaystyle{\lim_{n\rightarrow
\infty}}q^{nt}\left(
1+\frac{q^{-nk}}{[k]_qa^k}\right)_{q,k}^{\frac{t}{k}}\\\nonumber
&=&[k]_q^{\frac{t}{k}}a^t\left(1+\frac{1}{[k]_qa^k}\right)_{q,k}^{\frac{t}{k}}
\displaystyle{\lim_{n\rightarrow \infty}} q^{nt}\frac{\left(
\frac{1}{[k]_qa^k}+q^k \right)_{q,k}^{n}}{\left(
\frac{q^t}{[k]_qa^k}+q^k\right)_{q,k}^{n}}\\\nonumber
&=&[k]_q^{\frac{t}{k}}a^t\left(1+\frac{1}{[k]_qa^k}\right)_{q,k}^{\frac{t}{k}}
\displaystyle{\lim_{n\rightarrow\infty}}
\frac{(1+[k]_qa^kq^k)_{q,k}^{n}}{(1+[k]_qa^k q^{k-t})_{q,k}^{n}}
\\\nonumber
&=&[k]_q^{\frac{t}{k}}\frac{a^t}{1+[k]_qa^{k}}\left( 1+
\frac{1}{[k]_qa^k}\right)_{q,k}^{\frac{t}{k}} \left( 1+[k]_qa^k
\right)_{q,k}^{1-\frac{t}{k}}.\nonumber
\end{eqnarray}

Thus, $B_{q,k}(t,nk)=c(a,t)\beta_{q,k}^{(a)}(t,nk),$ for all
$n\in\mathbb{Z^+}.$ Moreover, proceeding as in Theorem
\ref{tildes} we can prove that part 2 of Theorem \ref{teor}, i.e.,
one can show that both sides of equation
\begin{equation*}
\displaystyle{B_{q,k}(t,s)=c(a,t)[k]_q^{-\frac{t}{k}}
\int_{0}^{\infty /a} \frac{x^{t-1}}{\left(1+
\frac{x^k}{[k]_q}\right)_{q,k}^{\frac{t+s}{k}}} d_qx, \quad
t,s>0,}
\end{equation*} after the appropriated changes, are formal power series with
rational coefficients in the correct variables. Part 1 follows
from part 2 using properties
$\beta_{q,k}^{(a)}(t,\infty)=(1-q)^{\frac{t}{k}}\gamma_{q,k}^{(a)}(t)$
and\\
$B_{q,k}^{(a)}(t,\infty)=(1-q)^{\frac{t}{k}}\Gamma_{q,k}^{(a)}(t).$
\end{proof}
\end{thm}
Below we include an alternative prove of this Theorem \ref{teor}.
Let us first give a proposition with further properties of the
function $c(a,t).$
\begin{prop}
\begin{enumerate}
\item $\displaystyle{\lim_{q\rightarrow1}c(a;t)=1}$ for all $a>0$
and $t
\in\mathbb{R}.$\\
\item $\displaystyle{\lim_{q\rightarrow0}c(a,t)=a^t+a^{t-k}}$ for
all $a>0$ and $0<t<1.$ \\
\item c(a,t) satisfies the following recursive formula: $c(a,t+k)=q^tc(a,t),$ for all $a>0$ and $t\in\mathbb{R}.$\\
\item For $a>0$ and $n\in\mathbb{Z^+},$ we have that $c(a,nk)=q^{kn(n-1)/2}.$\\
\item $\partial_qc(a,t)=0$, for all $a>0$ and $t\in \mathbb{R}.$
\end{enumerate}

\begin{proof}
\begin{enumerate}
\item Obvious.

\item In the limit $q\rightarrow0$, c(a,t) goes to
$\frac{a^t}{1+a^k}\left(
1+\frac{1}{a^k}\right)(1+a^k)=a^t+a^{t-k}$, for all  $a>0$ and $t
\in\mathbb{R}.$\\

\item $\:$ \vspace{-1cm}
\begin{align*}\frac{c(a,t+k)}{c(a,t)}&=[k]_q\frac{a^k\left(
1+\frac{q^t}{[k]_qa^k}\right)}{(1+[k]_qq^{-t}a^k)}=q^t.
\hspace{12cm}
\end{align*}
\item Immediate from item 3 and the fact that $c(a,0)=c(a,k)=1.$
\item To show  that $\partial_qc(a,t)=0,$ it is enough to check
that c(qa,t)=c(a,t).

\begin{align*}
c(qa,t)&=\frac{[k]_q^{\frac{t}{k}}q^ta^t}{1+[k]_qq^ka^k}\left(
1+\frac{1}{[k]_qq^ka^k}\right)_{q,k}^{\frac{t}{k}}(1+[k]_qq^ka^k)_{q,k}^{1-\frac{t}{k}}\\
&=\frac{[k]_q^{\frac{t}{k}}q^ta^t\left(
1+\frac{1}{[k]_qa^k}\right)_{q,k}^{\frac{t}{k}}\left(\frac{1}{[k]_qq^ka^k}+1
\right)(1+[k]_qa^k)_{q,k}^{1-\frac{t}{k}}(1+[k]_qq^{k-t}a^k)}{(1+[k]_qa^k)(1+[k]_qq^ka^k)\left(\frac{q^t}{[k]_qq^ka^k}+1\right)}.
\end{align*}

It is easy to check that $\displaystyle
{q^t\frac{\left(\frac{1}{[k]_qq^ka^k}+1
\right)(1+[k]_qq^{k-t}a^k)}{(1+[k]_qq^ka^k)\left(\frac{q^t}{[k]_qq^ka^k}+1\right)}=1},$
concluding that \ \vspace{0.2cm} c(qa,t)=c(a,t). Moreover
$\partial_qc(a,t)=0$, for all $t\in \mathbb{R}.$
\end{enumerate}
\end{proof}
\end{prop}

Theorem \ref{teor} may also be deduced from the following chain of
arguments. First, notice that using the Jackson integral,
Definition \ref{sepa} item 2 and the infinite product expression
for the function $B_{q,k}$ given in Lema \ref{betica} part 2, one
can show that Theorem \ref{teor} part 2, that is,
\begin{equation*}
\displaystyle{B_{q,k}(t,s)=c(a,t)[k]_q^{-\frac{t}{k}}
\int_{0}^{\infty /a} \frac{x^{t-1}}{\left(1+
\frac{x^k}{[k]_q}\right)_{q,k}^{\frac{t+s}{k}}} d_qx, \quad
t,s>0,}
\end{equation*}
is equivalent to the following relation
\begin{equation*}\label{sin}
\sum_{n \in \mathbb{Z}}\frac{q^{nt}\left(
1+\frac{1}{[k]_qa^k}\right)_{q,k}^{n}}{\left(1+\frac{q^{t+s}}{[k]_qa^k}\right)_{q,k}^{n}}=\frac{(1-q^k)_{q,k}^{\infty}(1-q^{t+s})_{q,k}^{\infty}\left(1+\frac{q^t}{[k]_qa^k}\right)_{q,k}^{\infty}\left(1+\frac{[k]_qq^ka^k}{q^t}\right)_{q,k}^{\infty}}{\left(1+\frac{q^{t+s}}{[k]_qa^k}\right)_{q,k}^{\infty}(1+[k]_qq^ka^k)_{q,k}^{\infty}(1-q^s)_{q,k}^{\infty}(1-q^t)_{q,k}^{\infty}}.
\end{equation*}
Making the changes  $u=-1/[k]_qa^k, v=-\frac{q^{t+s}}{[k]_qa^k}$
and $x=q^t$, we see that Theorem \ref{teor} part 2 is equivalent
to the famous Ramanujan identity.

\begin{equation*}\label{rama}
\sum_{n\in\mathbb{Z}}\frac{x^n(1-u)_{q,k}^{n}}{(1-v)_{q,k}^n}=\frac{(1-q^k)_{q,k}^{\infty}(1-v/u)_{q,k}^{\infty}(1-ux)_{q,k}^{\infty}\left(
1-q^k/ux\right)_{q,k}^{\infty}}{(1-v)_{q,k}^{\infty}(1-q^k/u)_{q,k}^{\infty}(1-x)_{q,k}^{\infty}(1-v/ux)_{q,k}^{\infty}}.
\end{equation*}

Similarly, using the definition of Jackson integral, Definition
\ref{sepa} part 2, the  infinite product expression for
$e_{q,k}^{-\frac{x^k}{[k]_q}}$ and Lemma \ref{ginfini}, one can
show that Theorem \ref{teor} part 1, that is,
\begin{equation*}
\displaystyle{\Gamma_{q,k}(t)=c(a,t)\int_{0}^{\infty/a(1-q^k)^{\frac{1}{k}}}
x^{t-1} e_{q,k}^{-\frac{x^k}{[k]_q}}}d_qx, \quad t>0,
\end{equation*}
 is equivalent to the following triple
identity
\begin{equation}\label{equ} (1-q^k)_{q,k}^{\infty}\left(
1+\frac{q^t}{[k]_qa^k}\right)_{q,k}^{\infty}\left(1+\frac{[k]_qq^ka^k}{q^t}\right)_{q,k}^{\infty}=(1-q^t)_{q,k}^{\infty}(1+[k]_qq^ka^k)_{q,k}^{\infty}\sum_{n\in\mathbb{Z}}q^nt\left(1+\frac{1}{[k]_qa^k}\right)_{q,k}^n.
\end{equation}
Making the change $x=-\frac{q^t}{a^k}$ and letting $a\rightarrow0$
in (\ref{equ}), we obtain
\begin{equation*}\label{jacobi} \sum_{n \in
\mathbb{Z}}(-1)^n
q^{kn(n-1)/2}[k]_q^{-n(n-1)/2}x^n=(1-q^k)_{q,k}^{\infty}\left(1-\frac{x}{[k]_q}\right)_{q,k}^{\infty}\left(
1-\frac{q^k[k]_q}{x}\right)_{q,k}^{\infty}
\end{equation*}  which is the Jacobi triple product identity.

\subsection* {Acknowledgment}
We thank Eddy Pariguan who has helped us in many ways.

$$\begin{array}{c}
  \mbox{Rafael D\'\i az. Instituto Venezolano de Investigaciones Cient\'\i ficas (IVIC).} \ \  \mbox{\texttt{radiaz@ivic.ve}} \\
  \mbox{Carolina Teruel. Universidad Central de Venezuela (UCV).} \ \  \mbox{\texttt{cteruel@euler.ciens.ucv.ve}} \\
\end{array}$$

\end{document}